\documentclass[12pt]{article}		

\def\preprint{TUW--02--10}
\def\finished{April 2002}
\def\archive {math.SC/0204356}

\long\def\Abstract{
	We describe our package PALP of C programs for calculations with 
	lattice polytopes and applications to toric geometry, which is freely
	available on the internet.
	It contains routines for vertex and facet enumeration,
	computation of incidences and symmetries, as well as completion of
	the set of lattice points in the convex hull of a
	given set of points. In addition, there are
	procedures specialised to reflexive polytopes such as the
	enumeration of reflexive subpolytopes,
	and applications to toric geometry and string theory, 
	like the computation of Hodge data and fibration structures
	for toric Calabi-Yau varieties.
	The package is well tested and optimised in speed as it was
	used 
	for time consuming tasks such as the classification of reflexive
	polyhedra in 4 dimensions and the creation and manipulation of
	very large lists of 5-dimensional polyhedra.
	While originally intended for low-dimensional applications, the
	algorithms work in any dimension and our key routine for vertex
	and facet enumeration compares well with existing packages.
}
\def\poly{{\it poly.x}}    \def\class{{\it class.x}} 	\def\nef{{\it nef.x}}
\def\glo{{\it Global.h}}				\def\cws{{\it cws.x}}
\def\PDm{{\tt POLY$_-$Dmax}}

\def\BI#1\EI{\setbox0\hbox{\tt #1}\unitlength=.2pt\BP(0,0)
	\put(-1,-1){\usebox0}\put(1,1){\usebox0}
	\put(-1,1){\usebox0}\put(1,-1){\usebox0}\EP}
\def\BI#1\EI{{\bf#1}}
\def\BL{\begin{itemize}\item[]\small\begin{alltt}}
\def\EL{\end{alltt}\end{itemize}}
\def\TO{$>$}  \def\PL{$|$}  \def\DQ{$\,'\hspace{-.5pt}'\,$}

\usepackage{cite,isolatin1,amssymb,latexsym}	\def\bye{\end{document}}
\long\def\new#1\endnew{{\bf #1}}  		\long\def\del#1\enddel{}
\def\ifundefined#1{\expandafter\ifx\csname#1\endcsname\relax}

\def\BC{\begin{center}} \def\EC{\end{center}} 	\catcode`\"=\active \let"=\"
\usepackage{alltt}

\let\0=\over      \let\1=\vec      \def\2{{1\over2}}
\let\4=\underline \let\5=\overline \let\6=\partial      \def\7#1{{#1}\llap{/}}
\def\8#1{{\textstyle{#1}}}         \def\9#1{{\ifmmode{\pmb{#1}}\else\bf#1\fi}}

\def\EEL#1 {\label{#1}\EE}		
\def\BE {\begin{equation}} 	\def\EE {\end{equation}}
\def\BEA{\begin{eqnarray}}	\def\EEA{\end{eqnarray}}

\let\and=\wedge                   \let\ex=\times

\let\bra=\langle        \let\ket=\rangle        \def\<#1\>{\bra #1 \ket}
	\def\rel#1 #2{\buildrel #1 \over {#2}}

\def\fnote#1#2{\begingroup\def\thefootnote{#1}\footnote{#2}
                \addtocounter{footnote}{-1}\endgroup}

     \let\o=\omega            
               
             \let\D=\Delta

\def\IR{{\mathbb R}}  \def\IP{{\mathbb P}}
  \def\IN{{\mathbb N}}
\def\IZ{{\mathbb Z}}

\def\plb#1 #2 {Phys. Lett. {\bf B#1} #2 }
\def\phr#1 #2 {Phys. Rep. {\bf  #1} #2 }        
\def\npb#1 #2 {Nucl. Phys. {\bf B#1} #2 }
\def\aph#1 #2 {Ann. Phys. {\bf #1} #2 }         
\def\jmp#1 #2 {J. Math. Phys. {\bf #1} #2 }
\def\jgp#1 #2 {J. Geom. Phys. {\bf #1} #2 }
\def\prd#1 #2 {Phys. Rev. {\bf D#1} #2 }
\def\prl#1 #2 {Phys. Rev. Lett. {\bf #1} #2 }
\def\rmp#1 #2 {Rev. Mod. Phys.  {\bf #1} #2 }
\def\zpc#1 {Z. Phys. {\bf #1C} }
\def\cmp#1 #2 {Commun. Math. Phys. {\bf #1} #2 }
\def\cqg#1 #2 {Class.Quant.Grav. {\bf #1} #2 }
\def\mpl#1 {Mod. Phys. Lett. {\bf A#1} }
\def\cpc#1 {Computer Phys. Commun. {\bf #1} }   
\def\ijmp#1 {Int. J. Mod. Phys. {\bf A#1} }
\def\ijmpC#1 {Int. J. Mod. Phys. {\bf C#1} }
\def\atmp#1 #2 {Adv. Theor. Math. Phys. {\bf #1} #2 }

\def\BP{\begin{picture}} \def\EP{\end{picture}}         
\newcounter{TRefNX} \let\OLDcite=\cite  \makeatletter
\def\makeTRefs#1{\@for  \NewTRef:=#1\do{\global\makeTRef{\NewTRef}}}
\def\makeTRef#1{\ifundefined{TRef#1}\stepcounter{TRefNX}%
\expandafter\xdef\csname TRef#1\endcsname{\theTRefNX}\fi}\makeatother
\def\NEWcite#1{\makeTRefs{#1}\OLDcite{#1}}

\ifundefined{draftmode} {} \else        \let\cite=\NEWcite
   \newcount\HOUR  \HOUR=\the\time \divide\HOUR by 60    \multiply \HOUR by 60
   \newcount\MIN   \MIN=\the\time  \advance\MIN by -\HOUR   \divide\HOUR by 60
   \ifnum   \MIN>9  \def\printTIME{{\it\the\HOUR\,:\,\the\MIN}}
         \else   \def\printTIME{{\it\the\HOUR\,:\,0\the\MIN}} \fi 
   
   \def\LLab#1{\BP(0,0)\unitlength=1mm\put(-12,.5){\makebox(0,0)[cr]{\small #1
        \rlap{$_{_{\makeatletter\csname TRef#1\endcsname\makeatother}}$}}}\EP}
\fi

\begin{document}

\vspace*{-18pt}\begin{flushright}  \archive\\[3pt] \preprint  \end{flushright}
\vspace{-3mm}

\BC
{\bf \Huge PALP:\\[9pt]
     \LARGE A Package for Analyzing Lattice Polytopes \\[9pt]
       with Applications to Toric Geometry
}
\\[9mm]
        Maximilian KREUZER\fnote{*}{e-mail: kreuzer@hep.itp.tuwien.ac.at}
\\[5mm]
        Institut f\"ur Theoretische Physik, Technische Universit\"at Wien\\
        Wiedner Hauptstra\ss e 8--10, A-1040 Wien, AUSTRIA
\\[5mm]                       and
\\[5mm] Harald SKARKE\fnote{\#}{e-mail: skarke@maths.ox.ac.uk}
\\[3mm] Mathematical Institute, University of Oxford\\
	24-29 St. Giles', Oxford OX1 3LB, ENGLAND
\\[7mm]
\vfil
{\bf ABSTRACT } \\[7mm]  \parbox{15cm}
{\baselineskip=14.5pt ~~~\Abstract} \EC
\vfil \noindent \preprint \\ \finished 
\setcounter{page}{0} \thispagestyle{empty} \newpage \pagestyle{plain}

\section{Introduction}

PALP is a package that contains a number of procedures for computing various
data on lattice polytopes.
It is freely available from our website \cite{KScy} under the GNU
licence terms.
It is written in C and provides the general purpose
program \poly, as well as a number of further main programs for special
tasks (these include \class, which was created for classifying
reflexive polyhedra, and \nef, which analyzes complete intersection
Calabi-Yau manifolds; both are useful in the context of toric
geometry).
While these programs cover many applications relevant for
someone working with lattice polytopes or toric Calabi--Yau spaces,
a variety of basic
routines from the package can be used to create new main programs for other
purposes.

\del
The most important of these routines will be described below,
while more details can be found in the on-line documentation \cite{KScy}.

The input for the program consists of the dimension and the number of points,
followed by an integer matrix of coordinates. Alternatively, a more abstract
(but often very efficient) format can be used, where the polytope is given in
terms of the intersection of positive halfspaces with the set of solutions
of linear diophantine equations, which we call {\it combined weight systems}
(CWS). Both formats will be explained in more detail at the beginning of
the next section.

The basic task of the package is the enumeration of vertices and
facet equations of the convex hull of a given set of points.
a fairly efficient completion of the lattice points
in the convex hull of a given set of points.
More sophisticated tasks are the computation of the incidence structure
of the face lattice and its dual, as well as the analysis of symmetries of
a polyhedron and the computation of a normal form.
Concerning symmetries it is
important to note that we work in a linear setting
(i.e., the lattice origin is always treated as a special point)
so that symmetries of a
polytope refer to $GL(\IZ)$ lattice automorphisms that leave a polytope
invariant. This also applies to
the normal form, which can be used to identify polytopes that are equal up
to a {\it linear} lattice automorphism. Nevertheless,
affine symmetries and equivalences can be analyzed by adding an auxiliary
dimension and putting the original polyhedron at distance one from the origin
(which has to be added to the polytope in order that the convex hull spans
the dimension).
\enddel

Basic tasks of the package include the enumeration of vertices and
facet equations of the convex hull of a given set of points,
an efficient completion of the set of lattice points
in this convex hull, and computation of the incidence structure
of the face lattice.
It is also possible to analyze symmetries of a polytope in terms of
$GL(\IZ)$ lattice automorphisms that leave it invariant; here it is
important to note that we work in a linear setting
(i.e., the lattice origin is a distinguished point).
In the same setting the computation of a normal form is provided;
this is a specific ordering of the vertices and choice of coordinate
system such that the normal forms of any two polytopes related by a
lattice automorphism are the same.

PALP places special emphasis on reflexive polytopes (i.e., lattice
polytopes whose duals are again lattice polytopes) and contains (or
rather, has evolved from) routines for their classification.
Reflexive polytopes play an important role in algebraic geometry where
they encode the combinatorial data of Calabi--Yau hypersurfaces in
toric varieties, leading to applications also in string theory and
conformal field theory.
PALP contains a number of procedures for applications in these
contexts, like the evaluation of Hodge data or the search for
reflexive sections of reflexive polytopes (corresponding to fibration
structures of the Calabi--Yau spaces).
We have tried, wherever possible, to keep the description in this
paper at the level of lattice polytopes such that background knowledge
in toric geometry and/or confomal field theory is required only for
understanding parts of section 4.

Our original aim in constructing the package was optimal performance in
terms of calculation time, memory efficiency and numerical stability
for reflexive polytopes in four dimensions.
We also used it to produce a number of new results on 5- and
6-dimensional reflexive polyhedra.
Testing our vertex and facet enumeration routines on higher (e.g. 12,
15, 20) dimensional polytopes revealed that they compete well with
those of existing packages \cite{cdd,qhull,porta,lrs}.
Like ``cdd'',  PALP is based on a double description
type of algorithm \cite{chalg} and uses a simple pivoting strategy which is
turned on for large dimension with a parameter. Similarly,
progress in the polyhedron analysis is monitored only in high dimensions
(for both parameters the defaults are set in {\it Vertex.c}).

Turning to the limitations of the package, we first have to point out that we
work with fixed precision because speed was one of our main concerns when
designing the package. We are using two data types, ``Long'' and ``LLong'',
which are usually set to ``long'' (32 bit) and ``long long'' (64 bit) in the
file {\it Global.h}.
``LLong'' is routinely used in numerically critical circumstances.%
\footnote{
	The most notable example is, perhaps, the function
	``EEV\_To\_Equation'', which computes a new facet equation by
	joining the intersection of two equations with a newly found
	vertex. Here it seems to be unavoidable that intermediate expressions
	can be of the order of squares of both input and output numbers.
}
In many examples with small coordinates (also in high dimensions) one
may try so save CPU time by using only 32 bit for ``LLong''. More
importantly, the analysis of very
large reflexive polytopes typically 
requires 64 bit in all operations in dimension 5 or higher. In any case, 
we implemented many tests that should alert the user in almost all cases 
where an overflow invalidates the result of a calculation.

Maximal dimensions, numbers of points, vertices and faces, etc., also have to
be fixed at compilation time. This is so because the whole package was
originally designed for polytopes of low dimensions (typically 4 to 6),
where reflexivity puts a priori limits to these numbers.
While facet enumeration still works well in higher dimensions,
some parts of the package, like the computation of incidences,
may run into problems with stack size, because the number of
faces may soon get very large. 

The package in its present status assumes the input polyhedron
to span the dimension. This, and many other limitations, could 
be avoided by combining routines that already exist in the package. Already
now, the user can obtain a set of equations that define the affine subspace
that is spanned by the set of input points, which should be sufficient to
proceed in many applications.
While PALP contains powerful routines for dealing with very large
numbers of polytopes (sorting, compression, data base structures etc.)
these applications are currently restricted to the reflexive case.
We hope that the present form of the package is only the starting
point and that many other applications will be created.
We will be happy to provide some help in this process and/or to
include user contributions in an appropriate way.

The structure of this paper is such that it should reflect growing
specialisation from general aspects of lattice polytopes via
reflexive polytopes to applications in algebraic
geometry and conformal field theory.
For many users it may be sufficient to look at section 2 where
we describe basic I/O formats and options of the main program \poly.
The formats are
designed such that the main programs can be combined in pipelines for
iterating certain operations. How this works is illustrated in section 3,
where we present two examples: the search for minimal reflexive polyhedra
that contain a given polytope, and the classification of reflexive
polyhedra in 3 dimensions.
Section 4 is mainly intended for readers who are interested in applications
to algebraic geometry and conformal field theory. Here we provide some
basics of toric geometry and discuss some options that can be used to compute
Hodge data and Poincar\'e polynomials related to transversal polynomials.
In section 5 we explain the structure of the package, as it is reflected
in the dependences given in the makefiles
and the definitions and routines 
in the main header file {\it Global.h}. This should provide a
basic understanding of how the parameters have to be set in various contexts,
and should also enable more ambitious users to design their own programs for
applications that are not covered in the present version of PALP.

\section{IO formats and options for \poly}

After obtaining our package from our web site \cite{KScy} and
following the instructions given there for unpacking,
it is possible to use \poly\ more or less immediately, simply by typing
``make poly'' and calling \poly\ with its help option.
In the following, as in all of our examples, user input is marked
by bold face.
\newpage
\BL
\BI	\$ poly.x -h	\EI
This is ``poly.x'':  computing data of a polytope P
Usage:   poly.x [-<Option-string>] [in-file [out-file]]
Options (concatenate any number of them into <Option-string>):
  h  print this information            | n  do not complete polytope or
  f  use as filter                     |      calculate Hodge numbers
  g  general output:                   | i  incidence information
     P reflexive: numbers of (dual)    | s  check for span property
       points/vertices, Hodge numbers  |      (only if P from CWS)
     P not reflexive: numbers of       | I  check for IP property
       points, vertices, equations     | S  number of symmetries
  p  points of P                       | T  upper triangular form
  v  vertices of P                     | N  normal form
  e  equations of P/vertices of P-dual | t  traced normal form computation
  m  pairing matrix between vertices   | V  IP simplices among vertices of P*
       and equations                   | P  IP simplices among points of P*
  d  points of P-dual                  |      (with 1<=codim<=# when # is set)
       (only if P reflexive)           | Z  lattice quotients for IP simplices
  a  all of the above except h,f       | #  #=1,2,3  fibers spanned by IP
  l  LG-`Hodge numbers' from single    |      simplices with codim<=#
       weight input                    | ## ##=11,22,33,(12,23): all (fibered)
  r  ignore non-reflexive input        |      fibers with specified codim(s)
  D  dual polytope as input (ref only) |    when combined: ### = (##)#
Input:    degrees and weights `d1 w11 w12 ... d2 w21 w22 ...'
          or `d np' or `np d' (d=Dimension, np=#[points]) and
              (after newline) np*d coordinates
Output:   as specified by options
\EL
Thus, if we want to find out about the vertices and equations of the
hyperplanes bounding a polytope, we can use \poly\ and specify the
options ``-ve''.
If no input file is
specified the program asks for input, either in the CWS form of degrees and
weights (see below), or in terms of a matrix of points whose convex
hull determines the polytope.
The program first
reads the size of the matrix (number of lines and number of columns) and
interprets the smaller of the two as the dimension. The matrix can therefore
be given as a line of column vectors or as its transposed, a column of line
vectors:
%
\BL
\$ \BI	poly.x -ve		\EI
Degrees and weights  `d1 w11 w12 ... d2 w21 w22 ...'
  or `#lines #colums' (= `PolyDim #Points' or `#Points PolyDim'):
\BI	3 4		\EI
Type the 12 coordinates as dim=3 lines with #pts=4 colums:
\BI	1 0 0 0 \EI
\BI	0 1 0 0	\EI
\BI	0 0 1 0 \EI
3 4  Vertices of P
    1    0    0    0
    0    1    0    0
    0    0    1    0
4 3  Equations of P
  -1  -1  -1     1
   0   0   1     0
   0   1   0     0
   1   0   0     0
\EL	
The resulting facet equations are given as lines
``$a_1~\ldots~a_n~~c$'' normalized
such that the $a_i$ have no common divisor and the
inequalities $\vec a\cdot\vec x + c \ge 0$ are fulfilled for all points of $P$.
If the inhomogeneous terms $c$ are positive for all equations then the origin
is in the interior of $P$.
In that case we say that $P$ has the `IP property'.
For such an `IP polytope'
the dual polyhedron $P^*$ is, by
definition, the convex hull of its vertices $\vec a/c$.

Reflexivity of a lattice polytope $P$ is equivalent to $P^*$ being a
lattice polytope, i.e. to $c=1$ for all facet equations. In that situation
the lines $\vec a$ can be interpreted as vertices of $P^*$ and
\poly\ omits the numbers $c$, indicating that the resulting matrix can be
interpreted as vertices of the dual polytope (in the transposed input format).
As an example we consider
\BL
\$ \BI	poly.x -e		\EI
Degrees and weights  `d1 w11 w12 ... d2 w21 w22 ...'
  or `#lines #colums' (= `PolyDim #Points' or `#Points PolyDim'):
\BI	3 4		\EI
Type the 12 coordinates as dim=3 lines with #pts=4 colums:
\BI	 3 -1 -1 -1   \EI
\BI	-1  3 -1 -1   \EI
\BI	-1 -1  3 -1   \EI
4 3  Vertices of P-dual <-> Equations of P
   1   0   0
   0   1   0
   0   0   1
  -1  -1  -1
\EL	
While the omission of $c=1$ in the reflexive case may be confusing in the
beginning, the advantage is that the resulting output can directly be used as
input of a calculation for the dual polytope.

This example brings us to the second input format. A simple way to obtain
lattice polytopes is as the set of non-negative solutions $\vec X\in\IZ^N$
to a system of $k$ linear diophantine equations $\vec w_j\cdot\vec X=d_j$
with non-negative $(w_j)_i=w_{ji}$ and $\sum_jw_{ji}>0$ for all $i$.
If these equations are independent, the resulting polytope has
support in the sublattice of
codimension $k$ in $\IZ^N$ defined by the equations.
Restricting the allowed equations to
the case $d_j=\sum_iw_{ji}$ implies the existence of at most one interior
point, namely $\vec X$ with $X_i=1$, which is always a solution to
this system of equations.
By shifting this generic solution to the origin and working with
the points $x_l=\sum_iB_l{}^i(X_i-1)$ transformed to some basis of the $N-k$
dimensional lattice that supports the polyhedron, we obtain a polytope
whose interior point (if it has one) is at the origin.
This construction provides our alternative input format in terms of
lists of numbers $d_1~w_{11}~\ldots~w_{1N}~d_2~w_{21}~\ldots~w_{kN}$, which
we call (combined) weight systems (for $k>1$) or (C)WS.
\del
provide our alternative input format,
which by construction can only be used for certain polytopes with at most one
interior point.
\enddel
It can be shown that all reflexive polytopes are subpolytopes
of polytopes of this kind \cite{crp}. While it may require some exercise to
feel comfortable with this format, it is easy to see that our last example
is the simplex defined by the weight system ``$4~~1~1~1~1$''. An example of
a combined weight system is ``$2~1~1~0~0~~2~0~0~1~1$'', which describes a
square. More details can be found in ref. \cite{ams}.

We now give a description of all the options of \poly\ that is more
detailed than the one in the help screen.
Note that if no option is specified the program by default chooses
the generic option ``-g'' (not recommended for large dimensions; see below).

\begin{description}	\def\Item[#1]{\item[\hbox to 22 pt{#1\hss}]}
\Item[-h]
The help screen is displayed.
\Item[-f~~~]
The filter flag switches off the prompt for the input data also in
case where input is read from the terminal (i.e. when no input file is
given). This is useful for building pipelines.
\Item[-g~~~]
The following data is displayed.
First the numbers $\# p$ and $\# v$ of lattice points and vertices,
respectively, in the format ``M: $\# p$ $\# v$''.
The remaining output depends on whether $P$ is reflexive.
In this case, the numbers $\# d$ and of $\# e$ of dual lattice points
and vertices are displayed in the format ``N: $\# d$ $\# e$'',
followed by information on the Hodge numbers of the corresponding
Calabi--Yau manifold
(in the case of a three dimensional polytope corresponding to a K3
surface, where the Hodge numbers are determined anyway, information on
the Picard number and the `correction terms' is given instead).
Otherwise the number $\# e$ of facets is shown as ``F: $\# e$''.
Using this option implies the completion of the set of lattice
points in the convex hull (``points'' in the help screen always means
``lattice points'').
In the reflexive case it also leads to the completion of the dual
polytope and the computation of the complete incidence structure which
is required for the calculation of the Hodge numbers.
For large dimensions (say, beyond 9) these tasks may result in an
extremely long response time or in a crash of the program.
In such a case, therefore,
it is strongly advisable to use other options, e.g. -nve,
if informations on the number of lattice points or Hodge numbers are not
required.
\Item[-p~~~~]
The lattice points of the polytope are displayed.
\Item[-v~~~~]
The vertices of the polytope are displayed.
\Item[-e~~~~]
The equations of the hyperplanes bounding the polytope are displayed.
\Item[-m~~~~]
One gets the $n_v\times n_e$ matrix with entries
$\vec a_j\cdot\vec v_i + c_j$, $1\le i\le n_v$, $1\le j\le n_e$, where
$n_v$, $n_e$ are the numbers of vertices and equations, respectively.
The elements of this `pairing matrix'
represent the lattice distances between the respective vertices and
facets.
\Item[-d~~~~] If the polytope is reflexive the lattice points of the
dual polytope are displayed.
\Item[-a~~~~]
This is a shortcut for ``-gpvemd''; it can be combined with any other
options.
\Item[-l~~~~]
This option is relevant for applications in the context of Landau-Ginzburg
models and will be explained in section 4.
\Item[-r~~~~]
Any input that does not correspond to a reflexive polytope will be ignored.
This is useful for filtering out reflexive polytopes from a larger list and
saves calculation time if one is interested only in reflexive polytopes.
\Item[-D~~~~]
The input is regarded as the dual polytope.
As this makes sense only in the reflexive case there is an error
message (but no exit from the program) for non-reflexive input.
This option is useful, in particular, if one wants to have
control over the order of the points in the $N$ lattice.
\Item[-n~~~~]
The completion of the set of lattice points is suppressed.
Hence the Hodge numbers cannot be calculated and the output will look
like the one for non-reflexive polytopes even in the reflexive case.
In particular, if the input is not of the (C)WS type, the number of
points may be displayed wrongly.
For large dimensional polytopes this option saves a lot of calculation
time.
\Item[-i~~~~]
Information on the incidence structure is displayed.
\Item[-s~~~~]
This option refers to the following property of (combined) weight systems.
As we saw above, in the higher dimensional embedding defined by a
weight system, the polytope is bounded by the inequalities $X_i\ge 0$.
We say that the polytope has the span property if the pullbacks of
the equations $X_i=0$ to the subspace carrying the polytope are spanned by
vertices of the polytope, i.e. if these equations correspond to
facets.
With ``-s'' a message is given if the (combined) weight systems does
not have this property.
\Item[-I~~~~]
There is a message if the polytope does not have
the origin of the coordinate system in its interior.
\Item[-S~~~~]
The output contains the following two numbers.
The first is the number of lattice automorphisms (elements of
$GL(n,\IZ)$) that leave the polytope invariant; note that each such
automorphism acts as a permutation on the set of vertices.
The second one is the number of permutations of the set of vertices
that leaves the vertex pairing matrix invariant (after taking into
account the induced permutations of facets).
This number can also be interpreted as the number of automorphisms of
the lattice generated by the vertices; it may be larger than the
number of symmetries in the given lattice.
\Item[-T~~~~]
A coordinate change is performed that makes the matrix of
coordinates of the points specified in the input upper triangular,
with minimal entries above the diagonal.
This may be useful for representing the polytope in a specific lattice
basis consisting of points of $P$, or for finding the volume of a
specific cone (if the generators of the cone are the first input
points, the volume will be the product of the entries in the diagonal
after the transformation).
\Item[-N~~~~]
This option leads to the calculation of a normal form of the polytope,
i.e. a matrix containing the vertices in a specific order in a
particular coordinate system, such that this output is the same for
any two polytopes related by a lattice automorphism.
This is useful, for example, if two polytopes are suspected to be
isomorphic because they are isomorphic if and only if their normal
forms are identical.
\Item[-t~~~~]
The calculation of the normal form involves determining the
pairing matrix, a normal form for the pairing matrix, an analysis of
which symmetries leave this normal form invariant, a preferred
ordering of the vertices and a conversion of the resulting vertex
coordinate matrix to upper triangular form.
The results of these steps, which may provide further useful
information on the structure of $P$, are displayed.
\Item[-V~~~~]
The IP simplices whose vertices are also vertices of the dual
polytope are displayed (an IP simplex is a possibly lower dimensional
simplex with the lattice origin in its relative interior;
for more information see, e.g., \cite{ams}).
\Item[-P~~~~]
The IP simplices whose vertices are lattice points of the dual
polytope are displayed.
This option should only be used if it is fairly clear that the dual
polytope does not have too many lattice points.
\Item[-Z~~~~]
This option only has an effect if combined with V or P. Then
the program computes a quotient
action that defines the sublattice of the intersection with the N lattice
with the respective linear subspace that is spanned by the vertices of
an IP simplex.
\del
\item[-\#~~~~]
Search for fibration structures, i.e. for lower dimensional reflexive
subpolytops, by checking the intersections of $P*$ with all linear subspaces
spanned by IP simplices for reflexivity. This misses fibrations whose
corresponding subspaces are only spanned by a combination of two or more
IP simplices. $\#$ is an integer between 1 and 3 and the codimension of the
IP simplices is restricted to $1\le$codim$\le\#$ (in contrast to our usual
policy this option has a side effect on P if combined with that option).
\item[-\#\#~~]
Compute all fibration structures of the given type: For 11, 22 and 33 all
reflexive sections of codimension 1, 2 and 3, respectively, are constructed.
For 12 and 23 all reflexive subpolytopes of codimension 1 and 2 that
themselves contain a reflexive subpolytope with relative codimension 1 are
constructed. The output is a polytope in the N lattice whose choice of
bases and whose order of points reflects the fibration structure.
\enddel
\item[-1,~-2,~-3,~~-11,~-22,~-33,~~-12,~-23~~]
Information on fibration structures of the Calabi--Yau manifold
corresponding to a reflexive polytope is displayed.
These structures correspond to reflexive subpolytopes of the dual ($N$
lattice) polytope $P^*$ that are intersections of the dual polytope with a
linear subspace of $N_\IR$; see, e.g., \cite{k3,fft,ams}.
As this is a time consuming task and the desired output format may
vary, we have programmed two different versions.

If a single number $\#\in\{1,2,3\}$ is specified, the intersections of
$P^*$ with all linear subspaces spanned by IP simplices are checked
for reflexivity.
This misses fibrations whose corresponding subspaces are only spanned
by a combination of two or more IP simplices.
The codimension of the IP simplices is restricted to
$1\le$codim$\le\#$ (in contrast to our usual policy this option has a
side effect on P if combined with that option).

If two numbers are specified, all fibration structures of the given
type are computed.
For 11, 22 and 33 all reflexive sections of codimension 1, 2 and 3,
respectively, are constructed.
For 12 and 23 all reflexive subpolytopes of codimension 1 and 2 that
themselves contain a reflexive subpolytope with relative codimension 1 are
constructed. The output is a polytope in the N lattice whose choice of
bases and whose order of points reflects the fibration structure.
\end{description}

Finally we turn to the the question of what can be done if something
goes wrong in the application of \poly.
We have designed our package in such a way that it should exit with
an error message rather than crash or display wrong results.
The two main sources for possible problems are inappropriately set
parameters and numerical overflows.
The most important settings of parameters all occur at the beginning
of \glo\ which is described in more detail in section 5.

Here are some typical error messages.
If we want to analyze the Calabi--Yau sixfold that is a
hypersurface in $\IP^7$, the following may happen.
\BL
\BI\$ 8  1 1 1 1 1 1 1 1\EI
Please increase POLY_Dmax to at least 7
\EL
In this case one should edit \glo, setting
\BL
#define                 POLY_Dmax       7       /* max dim of polytope      */
\EL
and compile again.
Similarly the program may ask for changes of other basic parameters, all of
which are defined within the first 60 lines of \glo.

In many cases we have implemented checks with the help of the
``assert'' routine, leading to error messages such as the following.
\BL
poly.x: Vertex.c:572: int Finish_IP_Check(PolyPointList *, VertexNumList *,
EqList *, CEqList *, INCI *, INCI *): Assertion `_V->nv<32' failed.
Abort
\EL
In this case one should look up line 572 of {\it Vertex.c},
\BL
    assert(_V->nv<VERT_Nmax);
\EL
This means that the value of {\tt $_-$V->nv} has risen above the value 32
assigned to {\tt VERT$_-$Nmax} in \glo\ and that the value of
{\tt VERT$_-$Nmax} should be changed correspondingly.
An assertion failure that
does not refer to an inequality involving a parameter or an allocation
failure is very likely to point to a numerical overflow.

While the possible problems mentioned above are related to parameters
set to values that are too low, excessively high values may also lead
to problems.
This may result in slowing down the calculation (in particular,
calculation time depends very sensitively upon whether {\tt
VERT$_-$Nmax} is larger than 64) or in running out of memory which may
lead to messages such as {\tt Allocation failure in Find$_-$Equations}.
Excessively high parameter values may even cause an
overflow of the stack, which often only results in a {\tt Segmentation fault}.
On many systems stack overflow is indistinguishable from real program errors
and the precise location of the crash would have to be found with a debugger.
But first one can try to recompile with, say, a smaller {\tt POINT$_-$Nmax}.

\section{Other programs and examples}

In addition to \poly\, our package contains the further main programs
\class, \cws\ and \nef. 
In the following we give a brief description and examples related to \class\ 
and in subsection 3.2 we also mention all we need 
to know about \cws, which creates weight systems and combined weight systems.
{\it nef.x} (by Erwin Riegler) will be discussed briefly in the second
paragraph of section 4.

The program \class\ provides an implementation of the algorithm for the
classification of reflexive polyhedra that is described
in refs. \cite{ams,crp}
(except for the very first steps of the algorithm which are
implemented in \cws).
In contrast to \poly\ this requires the program to
work not on one polyhedron at a time but rather on a set of polyhedra of
the same dimension, typically related by inclusion (possibly on
sublattices). This task required the definition of new data and I/O
structures. In particular, because of the huge number of reflexive polytopes
in 4 dimensions, we had to implement some data compression and a
corresponding binary I/O format. Moreover, since the data exceeded 2GB even
in their compressed form, it was necessary to
devise a data base structure for
storing the data in a set of smaller files on the hard disk in a sorted
form, with the top of a binary tree kept in the program memory
to enable a fast search for individual polytopes. Fortunately, the user
does not need to know anything about this except that the relevant file
names have to be provided in combination with the respective flags.

Because of the potentially extremely large data that have to be processed
the problem of stack overflow is more serious for \class. While \poly\ can
easily cope with polytopes that have serveral million points, it does not
make sense to consider the set of subpolytopes of such a huge polytope, and
hence \class\ should be compiled with a more moderate setting of the relevant
parameters.

With the help option of \class\ detailed information on the various items is
directly provided by the executable on typing the letter for the respective
option.
\BL
\BI	\$ class.x -h	\EI
This is  `class.x', a program for classifying reflexive polytopes
Usage:     class.x  [options] [ascii-input-file [ascii-output-file]]
Options:   -h          print this information
           -f or -     use as filter; otherwise parameters denote I/O files
           -m*         various types of minimality checks (* ... lvra)
           -p* NAME    specification of a binary I/O file (* ... ioas)
           -d* NAME    specification of a binary I/O database (DB) (* ... ios)
           -r          recover: file=po-file.aux, use same pi-file
           -o[#]       original lattice [omit up to # points] only
           -s*         subpolytopes on various sublattices (* ... vphm)
           -k          keep some of the vertices
           -c          check consistency of binary file or DB
           -M[M]       print missing mirrors to ascii-output
           -a[2b]      create binary file from ascii-input
           -b[2a]      ascii-output from binary file or DB
           -H*         applications related to Hodge number DBs (* ...cstfe)

Type one of [m,p,d,r,o,s,k,c,M,a,b,H] for help on options,
`g' for general help, `I' for general information on I/O or `e' to exit:
\EL
Rather than duplicating that information we therefore proceed with two
examples that illustrate how a combination of \poly\ and \class\ can be used
to accomplish certain goals. Further examples are provided in the next section
on applications, which however uses notions from toric geometry
and conformal field theory more freely.

\subsection{A search for maximal reflexive subpolytopes}

In our first example we consider the weight system ``7 1 1 1 1 1 2'' which
describes a non-reflexive polytope
and construct its largest reflexive subpolytope.%
\footnote
{	In geometrical terms this corresponds to a refexive blow-up of a
	weighted projective space, which in this case admits a transversal
	Calabi-Yau hypersurface whose Newton polytope is not reflexive.
	This can only occur in 5 or more dimensions.
}
We first put the weight system on a file and compute the generic data that
are obtained with the option ``-g'' (this is the default if no option is
specified):
\BL
\BI	\$ echo \DQ7 1 1 1 1 1 2\DQ \TO w.in	\EI
\BI	\$ poly.x w.in				\EI
7 1 1 1 1 1 2 M:496 10 F:7
\EL
The polyhedron thus has 496 points, 10 vertices and 7
facets. It
is not reflexive, because then the output line would contain the number of
dual points and vertices and the Hodge data instead of the ``F'' information.

Next we use the program \class\ to search for reflexive subpolytopes:
\BL
\BI	\$ class.x -o5 -po zw.5 w.in \EI
rec-dep<=5
0kR-0 0MB 0kIP 0kNF-0k 491_491 v14r10 f9r7 4b7 0s 0u 0n
7 1 1 1 1 1 2 R=1 +0sl hit=0 IP=5837 NF=1 (0)
Writing zw.5: 1+0sl 0m+0s 7b done: 0s
\EL
The option ``-o5'' requests that up to (recursion depth) 5 points are omitted
in searching for reflexive subpolytopes. The program tells us (some diagnostic
information and) that it found one reflexive subpolytope.
In order to complete its task  it had to check 5837
polytopes for the existence of an interior point.
The result is written in compressed binary format onto the file ``zw.5''.
In order to get it in readable form we need another run of \class:
\BL
\BI	\$ class.x -b2a -pi zw.5\EI
5 10
1 1 1 1 1 -1 -1 -1 -1 -1
0 3 0 0 0 -2 -2 5 -2 -2
0 0 3 0 0 -2 -2 -2 5 -2
0 0 0 3 0 -2 -2 -2 -2 5
0 0 0 0 3 5 -2 -2 -2 -2
np=1+0sl  5d  1v<=7 n<=7  1nv  1 0 1 7  0 0 0 0
\EL
(the last line is statistical information on the content of the binary file,
which we ignore). More information on this polytope can be obtained with
\poly:
\BL
\BI	\$ class.x -b2a -pi zw.5 \PL poly.x -fg \EI
M:491 10 N:8 7 H:2,0,450 [2760]
\EL
so indeed we found a reflexive polytope that has 5 points less.%
\footnote{
	The Hodge numbers of the corresponding Calabi--Yau 4-fold are
	$h_{11}=2$, $h_{12}=0$, $h_{13}=450$, and the Euler number is 2760.}
The structure of that polytope can be analyzed further by computing the
IP simplices, i.e. the linear relations among the vertices of the dual
polytope
\BL
\BI	\$ class.x -b2a -pi zw.5 \PL poly.x -fV \EI
5 7  vertices of P-dual and IP-simplices
    1   -1   -1   -1   -1    2   -1
    0    0    0    1    0   -1    0
    0    0    1    0    0   -1    0
    0    1    0    0    0   -1    0
    0    0    0    0    1   -1    0
-----------------------------------   #IP-simp=2
    2    1    1    1    1    1    0    7=d  codim=0
    1    0    0    0    0    0    1    2=d  codim=4
\EL
We thus recover the original weight system and an additional vertex that is
located opposite to the first vertex. By displaying the facet equations with
``poly.x -e'' one can see that the facet corresponding to this new
dual vertex is parallel to the facet at distance 2 that spoiled reflexivity
of our starting polyhedron, but now it is at distance 1. The IP simplices tell
us that the maximal reflexive subpolytope can actually be described by a
combined weight system:
\BL
\BI	\$ echo \DQ7 2 1 1 1 1 1 0  2 1 0 0 0 0 0 1\DQ \PL poly.x -fg \EI
7 2 1 1 1 1 1 0  2 1 0 0 0 0 0 1 M:491 10 N:8 7 H:2,0,450 [2760]
\EL
Remaining doubts about the equality of these polytopes could be eliminated
by comparing the normal forms that are computed by ``poly.x -N''.

\subsection{Classification of reflexive polytopes in three dimensions}

We now want to know how to use our programs for classifying all
reflexive polytopes in three dimensions \cite{c3d,ams}.
As a first step we have to find the relevant weight systems with 4
weights.
The algorithm for finding a finite set of candidates for such weight
systems was found in \cite{wtc} and is also described in \cite{ams}.
It is implemented in the program \cws, which creates weight systems with
the option ``-w'' and combined weight systems with the option ``-c'' (help
is available with ``cws.x -h''). 
Extra information on the resulting (C)WS can be obtained by using
\class\ with the option ``-ma'' ({\it m} for minimality, {\it a} for
all types of minimality checks), 
which is designed in such a way that it returns only those (C)WS that
have the IP property, together with information on the minimality type.
So the first steps are
\BL
\BI \$ cws.x -w3 \TO ws.3d \EI
\BI \$ cws.x -c3 \TO cws.3d \EI
\BI \$ cat cws.3d ws.3d \PL class.x -ma -f \TO wK3.ma \EI
\EL
\del
The file ``cws.3d'' now contains all 25 candidates for combined 
weight systems (21 of these have the IP property).
Next we also need the combined weight systems
After compiling {\it crecws}, we just type
\BL
\BI \$ cat c4 c3u3 \TO cK3 \EI
\EL
and adding the lines
\BL
3 1 1 1 0 0  2 0 0 0 1 1
4 2 1 1 0 0  2 0 0 0 1 1
6 3 2 1 0 0  2 0 0 0 1 1
2 1 1 0 0 0 0  2 0 0 1 1 0 0  2 0 0 0 0 1 1
\EL
corresponding to trivial combinations of lower-dimensional weight
systems (only non-trivial combinations are generated by {\it crecws}).
In order to find which of these candidates have the IP property and
certain minimality properties, 
Next we concatenate these files
and use \class\ with the option
``-ma'' ({\it m} for minimality, {\it a} for all types of minimality checks),
which is designed in such a way that it returns only those (C)WS that
have the IP property, together with information on the minimality type.%
\footnote{Mind the comment line with the numbers 99 of candidate weights 
	and 95 of resulting IP weights at the end of ``ws.3d''.
}
\BL
\BI \$ cat cws.3d ws.3d \TO cK3 \EI
\BI \$ class.x -ma cK3 \TO wK3.ma \EI
\EL
\enddel
This leads to the complete list of 116 IP (C)WS for $d=3$:
\BL
\BI \$ wc -l wK3.ma \EI
    116 wK3.ma
\EL
Every line in the resulting file contains extra information on the
minimality type which can be used to extract the 15 r-minimal (C)WS:
\BL
\BI \$ grep r wK3.ma \TO wK3.r \EI
\BI \$ wc -l wK3.r \EI
     15 wK3.r
\EL
Now we are ready to use \class\ for the construction of all subpolytopes
of these 15 polytopes:
\BL
\BI\$ class.x -po zK3 wK3.r \EI
...
Writing zK3: 4318+0sl 2119m+79s 11549b  u60 done: 0s
\EL
This means that \class\ has found 4318 subpolytopes on the original
lattices and no subpolytope that is reflexive only on a sublattice but
not on the original one.
The 4318 subpolytopes are encoded in the binary file {\it zK3}.
We can convert this file into a database,
\BL
\BI\$ class.x -pi zK3 -do zzK3 \EI
Read zK3 (4318poly +0sl 11549b)  write zzK3.* (11 files)  done (1s)
\BI\$ ls zz* \EI
zzK3.info  zzK3.v05   zzK3.v07   zzK3.v09   zzK3.v11   zzK3.v13
zzK3.v04   zzK3.v06   zzK3.v08   zzK3.v10   zzK3.v12
\EL
and look for polytopes in this database that are reflexive on a
sublattice such that the sublattice-polytope is not in the database.
\BL
\BI\$ class.x -sv -di zzK3 -po zzK3.s \EI
...
Writing zzK3.s: 1+0sl 0m+0s 3b done: 0s
\EL
There is precisely one such polytope (written to {\it zzK3.s}) and we
convert it to ascii-format by using the ``-b[2a]'' (``binary to
ascii'') option of \class:
\BL
\BI\$ class.x -b -pi zzK3.s \EI
3 4
1 1 1 -3
0 2 0 -2
0 0 4 -4
\EL
By computing the IP simplices among the vertices of the dual polytope
and the group action that defines the sublattice spanned by these vertices

\BL
\BI\$ class.x -b -pi zzK3.s \PL poly.x -fVZ \EI
3 4  vertices of P-dual and IP-simplices
    3   -1   -1   -1
   -2    2    0    0
   -1    0    1    0
--------------------   \#IP-simp=1
    1    1    1    1    4=d  codim=0 /Z2: 1 0 1 0
\EL
we see that the 4319{\it th} polytope is a $\IZ_2$ quotient of the
``quartic''. Our program, actually, accepts an extended version of the CWS
input, where a group action for selecting a sublattice can be added in the
format produced by the option Z:
\BL
\BI	\$ echo \DQ4 1 1 1 1 /Z2: 1 0 1 0\DQ \PL poly.x -fgN \EI
4 1 1 1 1 /Z2: 1 0 1 0 M:19 4 N:7 4 Pic:9 Cor:6
3 4  Normal form of vertices of P
   1   1   1  -3
   0   2   0  -2
   0   0   4  -4
\EL
The normal form indeed coincides with the one extracted from the file
zzK3.s above.

A complete list of all 4319 polytopes in ascii format can be generated
with the following steps:
\BL
\BI\$ class.x -pi zK3 -pa zzK3.s -po zK3.all \EI
Data on zK3:  4318+0sl  11825b  (3d)
Data on zzK3.s:  1+0sl  51b  (3d)
SL: 0nf 0sm 0nm 0b ->
d=3 v10 v<=13 n<=10 vn43  2199 79 0 11549  0 0 0 0
Writing zK3.all: 4319+0sl 2120m+79s 11549b   [p^2/2m=0M]
\BI\$ class.x -b -pi zK3.all \TO zK3.all.as \EI
\EL
With the data ``2120m+79s'' in the last output line mirror symmetry of the 
result is easily checked, as twice the number 2120 of mirror pairs and the 
number 79 of selfdual polytopes add up to 4319.

\section{Some applications in toric geometry and conformal field theory}

In the present section we present some further examples and results that
we produced while preparing the package for publication.
In particular we
analyzed the list of 5-dimensional Newton polytopes of the transversal
weight systems for Calabi--Yau 4-folds that were classified in \cite{lsw}.
This set contains polytopes with up to 355\,785
points and with degrees (and
coordinates) as large as 6\,521\,466, and thus provided stringent tests on
the numerical quality of our algorithms.%
\footnote{
        For these polytopes, single precision (i.e. ``Long'' set to 32
        bit integers) suffices
        until the
	square of the degree of the weight system approaches $10^{31}$. The
	only place where 64 bit are required much earlier is the routine
	``EEV\_To\_Equation'' in Vertex.c, which combines the intersection of
	two hyperlane equations and a newly found Vertex to a new hyperplane
	equation. This problem cannot be solved by an improvement of the
	algorithms and therefore some intermediate numbers in that routine
 	are always computed with maximal precision.
}
We discuss some new results on weight systems and what kind of
data for superconformal field theories (Landau Ginzburg models) can be
computed with \poly. Then we
give some examples of free quotients and fibration structures in toric
Calabi--Yau hypersurfaces.

The package PALP also contains procedures
that construct mirror pairs of toric complete intersections by enumerating
nef partitions of a reflexive polytope \cite{bb1,bb2} and
by evaluating the Hodge data of the corresponding manifolds
\cite{strh,nef,proc}. Since this
application, written by Erwin Riegler, is still in an experimental state,
we do not provide any documentation. It is, however, already quite reliable
and the information that is provided by the help option of the executable
\nef\ should be sufficient at least for the computation of some basic data.

\subsection{Landau--Ginzburg models and weighted projective spaces}

The first large list of Calabi--Yau manifolds with approximate mirror
symmetry consisted of hypersurfaces in weighted projective spaces
\cite{cls}. To avoid singularities that are not inherited from
the ambient space one requires transversality of the polynomial equations
that define the varieties.
Together with the Calabi--Yau condition $d=\sum \o_i$ this
makes the set of appropriate weight systems finite \cite{cqf,nms,KlS} .

The same quasi-homogeneus polynomials can be used as superpotentials of
Landau Ginzburg models \cite{lvw}, and it turns out that the Hodge data of
the manifold correspond to the charge degeneracies of the chiral ring of
the canonical LG orbifold \cite{va89,in90}. The basic object in this context
is the Poincar\'e polynomial \cite{lvw}
\BE
	P(t)=\prod {1-t^{d-\o_i}\01-t^{\o_i}},
\EE
in terms of which all orbifold data can be computed \cite{va89,in90,kr95}.

While it turned out that this class of models is incomplete because of its
lack of mirror symmetry (even if orbifolds and discrete torion are included
\cite{aas,dt}), Batyrev discoverd that the above construction can be
generalized to hypersurfaces in toric varieties  \cite{bat}.
He showed that the Calabi--Yau condition translates into reflexivity of the
polytopes, in terms of which the toric varieties are defined
\cite{Ful,cox,ck99}, and that the Hodge data can be computed in terms of the
incidence data and of the numbers of lattice points of dual pairs of faces
\cite{bat}. Starting from a weighted projective space, the corresponding
toric variety is found by considering the Newton polytope of the most general
quasi-homogeneous polynomial of degree $d=\sum\o_i$. That polytope is
automatically reflexive in up to 4 dimensions \cite{wtc}. (The simplest
counterexample in 5 dimension was the subject of our example in section 3.1.)

In a certain aspect, on the other hand, the Landau Ginzburg models are more
general. The resulting superconformal field theories can be used for
constructing consistent string vacua even if the condition of vanishing
canonical class is violated provided that the central charge
$c/3=\sum_i(1-2\o_i/d)$ is equal to the complex dimension of the ``internal
space''. Accordingly, the program \poly\ admits a different kind of input
when
called with the Landau--Ginzburg option ``-l''. While the input is restricted
to single weight systems, the degree need not be equal to the sum of the
weights (moreover, the degree may be given as the last instead of the first
number). If the resulting value $c/3$ for the central charge is an integer,
the analogue of the Hodge numbers is computed. Otherwise only the data of the
Newton polytope and the central charge are evaluated:
\hspace*{-12pt}\BL
\BI\$ echo \DQ3 1 1 1 1 1\DQ \PL poly.x -flg    \EI
1 1 1 1 1 3=d M:35 5 F:5 LG: c/3=5/3
\BI\$ echo \DQ3 1 1 1 1 1 1\DQ \PL poly.x -flg    \EI
1 1 1 1 1 1 3=d M:56 6 F:6 LG: H0:1,0,1 H1:0,20 H2:1
\EL
In the second case, the complete Hodge diamond is computed
but only part of it is displayed because the rest is fixed by
duality under the Hodge star operation; the output format for
$c/3=n\in \IN$ is
{\tt LG: H0:}$h_{00},\ldots,h_{0n}$
{\tt H1:}$h_{10},\ldots,h_{1,n-1}~~\ldots~~$ {\tt H$n$:}$h_{nn}$.

If the Newton polytope is reflexive, Vafa's formula for the CFT data is
compared with the geometrical data of the corresponding Calabi--Yau variety.
If there were a discrepancy the program would exit with a diagnostic message.
The fact that this never happend, even for the 252\,933 reflexive
among the 1\,100\,055 transversal weight systems for Calabi--Yau 4-folds, is
quite non-trivial. This coincidence, while natural, is not proven
nor really understood. From the computational point of view
Vafa's formula becomes quite slow at large degrees because the evaluation of
the Poincar\'e polynomial requires a division of polynomials whose degree
ranges up to 26 million in the case of 4-folds. Accordingly,
\BL
\BI
\$ echo \DQ6521466 1805 1806 151662 931638 2173822 3260733\DQ \PL poly.64 -lf
\EI
1805 1806 ... 3260733 6521466=d M:355 6 N:355785 6 V:303148,0,252 [1820448]
\EL
takes almost 10 minutes on an AMD XP1900 processor, while
\BL
\BI
\$ echo \DQ6521466 1805 1806 151662 931638 2173822 3260733\DQ \PL poly.64 -f\EI
6521466 1805 1806 ... 3260733 M:355 6 N:355785 6 H:303148,0,252 [1820448]
\EL
is finished in 8 seconds. With the ``-l'' flag the Hodge data are preceeded
by a ``V:'' to indicate that the numbers were computed with Vafa's formulas
(and compared with Batyrev's in the reflexive case, indicated by the ``N:'',
while a combination of ``F:'' and ``V:'' characterizes a nonreflexive
transversal weight system); moreover the degree of the weight system is
written after the weights because of the Landau--Ginzburg flag.

As one can verify with \class, distinct weighted projective spaces need not 
lead to different reflexive polytopes. With 
\BL
\BI\$ cws.x -w4 -t \TO tws.4d \EI
\EL
we create, as an example, a file containing all 7555 transversal weight 
systems for 4-dimensional Newton polyhedra (this takes a few hours, 
but the data can also be fetched from our web page \cite{KScy}). Then
\BL
\BI\$ class.x -a2b -po ztws.4d  tws.4d twsred.4d \EI
...
Writing ztws.4d: 5961+0sl 404m+32s 76243b  u32 done: 0s
\EL
writes the binary file {\it ztws.4d} encoding a sorted list of
5961 distinct polytopes and the ascii file {\it twsred.4d} containing a
reduced list of 5961 weight systems for these polytopes.
Among the 252,933 reflexive polytopes coming from
5-dimensional weighted projective spaces obeying the transversality
condition,
202,746 are different.
Similarly, the 184,026 `IP weight systems' with 5 weights correspond
only to 58,690 different 4-dimensional polytopes.

\subsection{Quotients and fibrations}

In toric geometry a refinement of the $N$ lattice, with the corresponding
change to a sublattice of the dual $M$ lattice, corresponds to the modding out
of an abelian group action. The resulting quotient of a Calabi--Yau
hypersurfaces is free, and thus leads to a non-trivial fundamental group, if
it does not introduce singularities. This is the case if the refinement of
the $N$ lattice does not lead to additional lattice points for the dual
polytope $\D^*\subset N_\IR$. With the option ``-sp'' the program \class\
searches either polytopes in ascii-format or the content of a data base for
such objects.%
\footnote{
	In fact a slightly weaker condition is sufficient: additional interior
	points of facets correspond to divisors that do not intersect the
	hypersurface. With the option ``-sh'' \class\ omits these points when
	checking for free quotients, but so far we did not find an example
	where this leads to additional solutions.
}
By searching the database of all reflexive 4-dimensional polytopes with
\BL\BI \$ class.x -sp -di zzdb	\EI\EL
we found that there are 16 free quotiens, for 13 of which the index of the
lattice quotient is 2. The remaining 3 quotients are the well-know free
$\IZ_5$ quotient of the quintic, a $\IZ_3$ quotient of $\IP^2\ex\IP^2$, and
a $\IZ_3$ quotient of weighted $\IP^4_{33111}$.

In string theory fibration structures play an important role in dualities.
In particular, K3 fibrations are essential for heterotic--type II duality
\cite{klm,AL,cafo} and elliptic fibrations are essential in F-theory
\cite{Fv}. In the context of toric geometry such fibrations correspond to
reflexive sections of codimension 1 or 2 of 4-dimensional polytopes in the
$N$ lattice \cite{k3,fft,hly}. In the remainder of this section we show
some examples of such fibrations that have a non-trivial fundamental group.%
\footnote{	It turned out that
	these quotient spaces admit an elliptic K3 fibration if and only 
	if the group is $\IZ_2$; this fact might find an explanation in the
	context of heterotic--type II duality \cite{kkrs}.
}
The only polytope in the list of 16 free quotients that does not have any
fibration is the $\IZ_5$ quotient. The two $\IZ_3$ quotients both are
elliptically fibered:
\BL
\BI\$ echo \DQ3 1 1 1 0 0 0  3 0 0 0 1 1 1 /Z3: 0 1 2 0 1 2\DQ\PL poly.x -fg2P\EI
3 1 1 1 0 0 0  3 0 0 0 1 1 1 /Z3: 0 1 2 0 1 2 M:34 9 N:7 6 H:2,29 [-54]
4 7  points of P-dual and IP-simplices
   -3    0    0    0    0    3    0
    1    1    0    0    0   -2    0
    1    0    0   -1    1   -1    0
    2    0    1   -1    0   -2    0
------------------------------    #IP-simp=2
    1    1    0    0    0    1   3=d  codim=2
    0    0    1    1    1    0   3=d  codim=2
------------------------------    #fibrations=2
    v    v    _    _    _    v  cd=2  m:10 3 n:4 3
    _    _    v    v    v    _  cd=2  m:10 3 n:4 3
\BI\$ echo \DQ9 3 3 1 1 1 /Z3: 1 2 1 2 0\DQ\PL poly.x -fg22PZ\EI
9 3 3 1 1 1 /Z3: 1 2 1 2 0 M:49 5 N:7 5 H:2,38 [-72]
4 7  points of P-dual and IP-simplices
   -3    0    3    0    0    0    0
    1    0   -2    0    3    1    0
    1    0   -1    1   -1    0    0
    2    1   -2    0   -1    0    0
------------------------------    #IP-simp=2
    3    1    3    1    1    0   9=d  codim=0 /Z3: 1 2 2 1 0 0
    1    0    1    0    0    1   3=d  codim=2
4 6  m:10 3 n:4 3  M:49 5 N:7 5  p=025134
   1   -2    1    0    0    3
   0    1   -1   -1    0   -2
   0    0    0    0    3   -3
   0    0    0   -1   -1    2
\EL
In the first case we used IP simplices to search for fibrations of
codim$\le2$. In 
the second case we can be sure that all elliptic
fibrations were found because of the 22
flag; the upper-left 2x3 block of the coordinate matrix describes
the boundary points of the fiber polytope in
the N lattice, whose dual has 10 points and 3 vertices (the information on
the fiber polytopes is given by ``m:10 3 n:4 3''). 
	The string ``p=025134'' at the end of the first line of the polytope
	in the fiber basis describes the permutation that brings the points 
	from the original order, to which the IP simplex information refers, 
	into an order where the points of the fiber come first; as usual in
	{\em C}, counting begins with 0. 

As our last example we consider a polytope with different elliptic K3
fibrations, where none of the K3 fibers is spanned by a single IP simplex:
\BL
\BI\$ poly.x -12PD	\EI
Degrees and weights  `d1 w11 w12 ... d2 w21 w22 ...'
  or `#lines #colums' (= `PolyDim #Points' or `#Points PolyDim'):
\BI4 7\EI
Type the 28 coordinates as dim=4 lines with #pts=7 colums:
\BI    1   -1    0    0    0    1   -1\EI
\BI    0    2   -1    0    2    1   -1\EI
\BI    0    0    0    1   -1    1   -1\EI
\BI    0    0    0    0    0    2   -2\EI
4 9  points of P-dual and IP-simplices
    1   -1    0    0    0    1   -1    0    0
    0    2   -1    0    2    1   -1    1    0
    0    0    0    1   -1    1   -1    0    0
    0    0    0    0    0    2   -2    0    0
----------------------------------------    #IP-simp=4
    1    1    2    0    0    0    0    0   4=d  codim=2
    0    0    2    1    1    0    0    0   4=d  codim=2
    0    0    1    0    0    0    0    1   2=d  codim=3
    0    0    0    0    0    1    1    0   2=d  codim=3
4 8  Em:9 3 n:5 3  Km:35 5 n:7 5  M:53 10 N:9 7  p=01273456
   1    1   -1    1    0    2    0    0
   0    2   -1    1    0    2    1   -1
   0    0    0    0    1   -1    1   -1
   0    0    0    0    0    0    2   -2
4 8  Em:9 3 n:5 3  Km:27 6 n:7 5  M:53 10 N:9 7  p=01275634
   1    1   -1    1    0    0    0    2
   0    2   -1    1    0    0    1    1
   0    0    0    0    1   -1    1   -1
   0    0    0    0    0    0    2   -2
...
4 8  Em:9 4 n:5 4  Km:27 6 n:7 5  M:53 10 N:9 7  p=25673401
   1    0    0   -1    0   -2    1   -3
   0    1   -1    0    0    0    1   -1
   0    0    0    0    1   -1    1   -1
   0    0    0    0    0    0    2   -2
\EL
The fibers can be checked to be elliptic curves in $\IP^1\ex\IP^1$
(for ``Em:9 4 n:5 4'') or in $\IP^2_{112}$, and the data of the elliptic
K3 surfaces are
``M:35 5 N:7 5 Pic:2 Cor:0'' or ``M:27 6 N:7 5 Pic:4 Cor:1'', as one can find
out by entering the upper-left 3x6 blocks of K3 coordinates into \poly\
with flags ``-gD''. 

The dots in the above result indicate 3 further fibrations, with 
corresponding permutations ``p=23470156'', ``p=23475601'', ``p=25670134'',
that can be checked to be related to the displayed fibrations  
by a symmetry of the polytope: The IP simplex structure is invariant under 
the four independent permutations (01), (34), (56) and (03)(14) of order 2, 
which induce the required relations. Indeed, ``poly.x -t'' displays the 
resulting 16 vertex permutations that represent all GLZ-symmetries of the 
polytope.
Note that in this case the double cover
\BL
\BI\$ poly.x -Dg	\EI
Degrees and weights  `d1 w11 w12 ... d2 w21 w22 ...'
  or `#lines #colums' (= `PolyDim #Points' or `#Points PolyDim'):
\BI4 7\EI
Type the 28 coordinates as dim=4 lines with #pts=7 colums:
\BI    1   -1    0    0    0    1   -1\EI
\BI    0    2   -1    0    2    1   -1\EI
\BI    0    0    0    1   -1    1   -1\EI
\BI    0    0    0    0    0    1   -1\EI
M:105 10 N:9 7 H:5,85 [-160]
\EL
has a larger Picard number $h_{11}=5$, which indicates that not all divisors
of the double cover can be toric. The number of points in the N lattice is
the same and the Euler number doubles, as it must be for a free quotient.
%
%
%
\del math.AG/0003033, Macaulay 2 and the geometry of schemes
\\class.x -a -po zz.aux TransRef.all TransRef.red
\\1 1 1 14 33 50 100=d M:23644 10 N:27 7 H:15,0,20155 [121068]
\\1 1 2 14 32 50 100=d M:12226 10 N:28 7 H:16,3,10439 [62760]
\\1 1 3 12 33 50 100=d M:9213 11 N:23 8 H:13,0,7855 [47256]
\\1 1 4 12 32 50 100=d M:7144 11 N:22 8 V:12,7,6099 [36672]
\\1 1 5 10 33 50 100=d M:6641 7 N:21 7 V:13,8,5661 [34044]
\\1 1 6  9 33 50 100=d M:6151 13 N:20 8 H:11,0,5245 [31584]
\\1 2 2 12 33 50 100=d M:6965 12 N:27 9 H:15,0,5973 [35976]
\\1 2 2 13 32 50 100=d M:6609 15 N:31 9 H:19,0,5647 [34044]
\\1 2 3 11 33 50 100=d M:5029 14 N:21 9 H:13,0,4286 [25842]
\\1 2 7  7 33 50 100=d M:3396 11 N:24 9 H:12,0,2896 [17496]
\\1 4 4  8 33 50 100=d M:2625 8 N:27 7 V:12,132,2240 [12768]
\enddel

\section{Structure of the package}

In addition to {\it Makefile} 
which controls the
compilation of executables, our package consists of program files
with extension {\it .c} and header files (extension {\it .h}).

\del
{\it GNUmakefile} is the default for GNU-make;
the options for {\tt gcc} or {\tt g++} used there should provide good 
performance.
\enddel
{\it Makefile} contains various different sets of options which are
optimized for different types of compilers 
(including GNU C++, Silicon Graphics and the DEC Alpha compilers).
For elaborate applications it may certainly be useful to experiment
with some of the options. 
Optimization level {\tt -O3} may result in faulty compilation
(we included several work-arounds for the GNU compiler
and for the Silicon Graphics), while {\tt -O2} is usually safe and only
marginally less effective.

The files in our package with extension {\it .c} come in two varieties:
the files starting with lower case letters (i.e. {\it class.c},
{\it cws.c}, {\it nef.c}, {\it poly.c}) correspond to main
programs whereas all other {\it .c}--files start with capitals and
contain routines used by the main programs.

{\it Coord.c}, {\it Polynf.c} and {\it Vertex.c} are used by 
all main programs.
Definitions and routines from these files that are globally available
are all listed in the header file {\it Global.h}, which is probably
the most interesting file for anyone trying to understand the
structure of our package or to do some programming with our routines.
As there is a lot of documentation within this file, we describe it
here only briefly and recommend also to study a printout of it.

{\it Global.h} starts with basic parameter settings.
This part is important not only for programmers but also for users:
none of the parameters should be chosen too small, but choosing them
too large may result in bad performance or even a crash due to lack of
memory.
Usually it should be sufficient to make sure that the
definition of \PDm\ (the maximal dimension of a polytope to be
analyzed) matches the application, because {\it Global.h} contains
default settings for other relevant parameters depending on \PDm.
For \PDm\ $\le$ 4 these settings are such that every
reflexive polytope fulfills them.
Next there are most of the type definitions that are used globally in our
package.
Finally, the globally available routines from {\it Coord.c}, {\it
Polynf.c} and {\it Vertex.c} are listed.

{\it Coord.c} contains routines mainly connected with I/O; these are
routines for reading coordinates (or CWS and converting these to
coordinates) as well as routines for various types of output.
Whenever possible, the output is structured in such a way that it can
be used as input for another run of our programs.

{\it Polynf.c} contains routines for analyzing the symmetries of a polytope,
calculating its normal form, identifying IP simplices, finding
fibrations, transforming bases and other advanced tasks.

{\it Vertex.c} contains the core routines for handling polytopes.
These include basic tasks such as evaluating an equation on a point,
comparing points, converting a codimension 2 subspace (defined by
two equations) and a vertex into the equation of the hyperplane
spanned by them, etc.
Next there are the routines of polytope analysis such as
determining the vertices and facets of the convex hull of a finite
number of lattice points and the computation of the complete set of
lattice points of a given polytope; there is also a routine for
calculating the Hodge numbers of a Calabi--Yau hypersurface determined
by a reflexive polytope.
Finally there are routines related to incidence structures.
In our programs these structures are represented by bit patterns which
are determined either by unsigned (long long) integers or by arrays of
unsigned integers.
This is useful because many of the bitwise operations provided by
C correspond to meaningful operations in terms of the analysis of a
polytope (see sec. 3 of \cite{c4d}).

{\it Subadd.c},  {\it Subdb.c} and {\it Subpoly.c} are only used by \class.
They share a common header file {\it Subpoly.h}.
{\it Subpoly.c} contains the routines necessary for constructing all
reflexive subpolytopes of a given polytope as well as many other
routines necessary for the classification of reflexive polytopes.
{\it Subadd.c} and {\it Subdb.c} both have to do with data handling.
In {\it Subdb.c} there are the routines concerning databases
consisting of binary files whereas  {\it Subadd.c} deals with other
tasks such as sorting lists of normal forms of polyhedra etc.

{\it LG.c} and the corresponding header file {\it LG.h} are used by
\poly\ and \cws, but not by \class.
{\it LG.c} contains the routines necessary for a CFT type analysis as
described in section 4.
\nef\ uses routines from {\it Nefpart.c} for creating nef partitions
and routines from {\it E\_Poly.c} for Hodge numbers; the corresponding
header file is {\it Nef.h}.

There is also a file {\it Rat.c} providing routines for dealing with
rational numbers and certain tasks concerning manipulations of integers.
The corresponding header file {\it Rat.h} is included wherever these
applications are needed.

$Acknowledgements.$ We would like to thank Victor Batyrev for drawing our
attention to free quotients and for suggesting the weak version of the
criterion that is implemented in \class. 
We also thank Erwin Riegler for permission to include his application
\nef\ in our package and for his contributions to \cws.
This work was supported in part by
the Austrian Research Funds FWF under grant Nr. P14639-TPH, and by the
EPSRC under grant reference GR/R61017/01 [HBKPF].


\bye